\documentclass[11pt]{amsart}

\newcommand {\R}{\mathbb{R}}

\newcommand {\grad}{\nabla}
\newcommand {\eps}{\epsilon}
\newcommand {\Dtilde}{\tilde{D}}
\newcommand {\Dbar} {\bar{D}}
\newcommand {\utilde}{\tilde{u}}
\newcommand {\abar} {\overline{\alpha}}
\newcommand {\calN} {\mathcal N}
\newcommand {\calS} {\mathcal S}
\newcommand {\weakto} {\rightharpoonup}

\newcommand {\dist}{{\rm dist\,}}
\newcommand {\supp}{{\rm supp\,}}

\newcommand {\pf}[1][Proof] {\vspace{\baselineskip}
                             \noindent {\em #1.} \quad}  

\newtheorem{theorem}{Theorem}
\newtheorem{prop}[theorem]{Proposition}

\newtheorem{cor}[theorem]{Corollary}
\newtheorem{lem}[theorem]{Lemma}
\newtheorem{conj}{Conjecture}
\newtheorem{prob}{Problem}

\begin{document}

\title[Symmetry breaking for composite membranes]
{Symmetry breaking and other phenomena in the optimization of eigenvalues for
composite membranes}

\author{S. Chanillo \and D. Grieser \and  M. Imai \and K. Kurata
      \and I. Ohnishi}
\renewcommand{\addresses}{\smallskip \noindent
Sagun Chanillo:
Department of Mathematics,
Rutgers University,
New Brunswick, NJ 08903, USA\\
chanillo@math.rutgers.edu
\vspace{1mm}

\noindent
Daniel Grieser:
Institut f\"ur Mathematik,
Humboldt-Universit\"{a}t Berlin,
Unter den Linden 6,
10099 Berlin, Germany\\
grieser@mathematik.hu-berlin.de
\vspace{1mm}

\noindent
Kazuhiro Kurata:
Department of Mathematics,
Tokyo Metropolitan University,
Minami-Ohsawa 1-1, Hachioji-shi, Tokyo, Japan\\
kurata@comp.metro-u.ac.jp

\noindent
Masaki Imai and Isamu Ohnishi:
Department of Infomation mathematics and Computer sciences,
University of Electro-Communications,
Chofu-ga-oka 1-5-1, Chofu-shi, Tokyo, Japan\\
imai-m@kenks.im.uec.ac.jp, ohnishi@im.uec.ac.jp
}

\date{February 20, 2000}

\begin{abstract}
We consider the following eigenvalue optimization problem:
Given a bounded domain $\Omega\subset\R^n$ and numbers $\alpha\geq 0$,
$A\in [0,|\Omega|]$, find a subset $D\subset\Omega$ of area $A$
for which the first Dirichlet eigenvalue of the operator
$-\Delta + \alpha \chi_D$ is as small as possible.

We prove existence of solutions and investigate their qualitative properties.
For example, we show that for some symmetric domains (thin annuli
and dumbbells with narrow handle) optimal solutions must possess fewer
symmetries than $\Omega$; on the other hand, for convex $\Omega$
reflection symmetries are preserved.

Also, we present numerical results and formulate some
conjectures suggested by them.
\end{abstract}

\keywords{
Optimal shape, minimization problem, first Dirichlet eigenvalue, 
free boundary, symmetry breaking, dumbbells, annulus}

\maketitle 

\section{Problem and Main Results} \label{secintro}
We study qualitative properties of solutions of
a certain eigenvalue  optimization problem. In physical terms,
the problem can be stated as follows:
\begin{quote}  {\bf Problem (P) }
Build a body of  prescribed shape out of given materials
(of varying densities) 
in such a way that the body has a prescribed
mass and
so that the basic frequency of the resulting
membrane (with fixed boundary) is as small as possible.
\end{quote}
In fact, we will consider a more general problem, which we now
state in mathematical terms:
Given a domain $\Omega\subset\R^n$ (bounded, connected, with Lipschitz boundary)
and  numbers $\alpha>0$, $A\in [0,|\Omega|]$ (with $|\cdot|$ denoting
volume).
For any measurable subset $D\subset\Omega$ let $\chi_D$ be its characteristic
function and $\lambda_\Omega(\alpha,D)$ the lowest eigenvalue $\lambda$ of the
problem
\begin{align} \label{eqeveq}
\begin{split}
-\Delta u + \alpha \chi_D u &= \lambda u \quad \text{ on }\Omega \\
u & = 0 \quad\text{ on } \partial \Omega.
\end{split}
\end{align}
Define
\begin{equation} \label{eqevmin}
 \Lambda_\Omega(\alpha,A) = \inf_{\substack{ D\subset\Omega \\
                                            |D| = A}}
                              \lambda_\Omega(\alpha,D).
\end{equation}
Any minimizer $D$ in \eqref{eqevmin} will be called an
{\em optimal configuration} for the data $(\Omega,\alpha,A)$.
If $D$ is an optimal configuration 
and $u$ satisfies \eqref{eqeveq} then $(u,D)$
will be called an {\em optimal pair} (or {\em solution}).
Our problem now reads:
\begin{quote} {\bf Problem (M)}
Study existence, uniqueness and qualitative properties of optimal pairs.
\end{quote}
 As is well-known,
$u$ is uniquely determined, up to a scalar multiple, by $D$,
and may be chosen to be positive on $\Omega$. In addition, we will
always assume $$\int_\Omega u^2 = 1.$$
(Integrals over $\Omega$ are always taken with respect to the standard measure.)
Clearly, changing $D$ by a set of measure zero does not affect 
$\lambda_\Omega(\alpha,D)$ or $u$. Therefore, we will consider sets $D$
that differ by a null-set as equal.

At first sight, it is not obvious that problem (M) generalizes
problem (P). 
In fact, we will see (Theorem \ref{thprobPM}) 
that there is a number $\abar_\Omega(A)>0$
such that solutions of problem (P) are in one to one correspondence
with solutions of problem (M) with parameters in the range
$\alpha\leq \abar_\Omega(A)$.
The number $\abar_\Omega(A)$ is characterized as the unique value of $\alpha$
satisfying
\begin{equation} \label{eqabardef0}
 \Lambda_\Omega(\abar_\Omega(A),A) = \abar_\Omega(A),
\end{equation}
see Proposition \ref{propparam}.
 
Our investigations are theoretical and numerical: Numerical results
(obtained by M.I.\ and I.O.) suggest properties of optimal configurations;
this leads to the formulation of conjectures, and some of these
are proved rigorously (by  S.C., D.G.\ and K.K.).

A central tool in our 
 investigations is the variational characterization
of the eigenvalue:
$$ \lambda_\Omega(\alpha,D) = \inf_{u\in H_0^1(\Omega)} R_\Omega(u,\alpha,D),
\qquad R_\Omega(u,\alpha,D) :=
     \frac{\int_\Omega |\grad u|^2 + \alpha\int_\Omega \chi_D u^2}{\int_\Omega u^2},$$
and the eigenfunction $u$ is a minimizer.
 So $\Lambda_\Omega(\alpha,A)$
is characterized by
$$ \Lambda_\Omega(\alpha,A) = \inf_{\substack{ 
                                     u\in H_0^1(\Omega) \\
                                     |D|=A}}
                                 R_\Omega(u,\alpha,D).   $$

We first prove the following theorem on existence and basic properties
of solutions. It is fundamental for all further considerations.

\begin{theorem} \label{thexist}
For any $\alpha>0$ and $A\in [0,|\Omega|]$ there exists an optimal pair.
Moreover, any optimal pair $(u,D)$ has the following properties:
\begin{enumerate}
\item[(a)] $u\in C^{1,\delta}(\Omega) \cap H^2(\Omega)
                \cap C^\gamma(\overline{\Omega})$
  for some $\gamma>0$ and every $\delta < 1$.
\item[(b)] $D$ is a sublevel set of $u$, i.e.\ 
there is a number $t\geq 0$ such that
$$D  = \{u\leq t\}.$$
\item[(c)] Every level set $\{u=s\}$, $s\geq0$, has measure zero, except
possibly in the case $\alpha=\abar_\Omega(A)$, $s=t$.
\end{enumerate}
\end{theorem}
Here we use the short notation $\{u=t\} = \{x: u(x) = t\}$.
Since $\chi_D$ is discontinuous, solutions $u$ may not be twice
differentiable, so equation \eqref{eqeveq} is understood in the weak sense.

Note that Theorem \ref{thexist}(b) 
shows in particular that our problem is equivalent
to finding the smallest eigenvalue and associated eigenfunctions
of the nonlinear problem (with free variables $u$ and $t$)
\begin{align} \label{eqnonlinearev}
\begin{split}
-\Delta u + \alpha \chi_{\{u\leq t\}} u &= \lambda u \quad \text{ on }\Omega \\
u & = 0 \quad\text{ on } \partial \Omega \\
|\{u\leq t\}| &= A.
\end{split}
\end{align}

The question of {\em uniqueness} is much more subtle: For some domains
$\Omega$ there will be a unique optimal pair for all $\alpha,A$,
while for others there will be many, for certain ranges of $\alpha,A$.
This follows from our
results on symmetry preservation and symmetry breaking below.

We now list a few questions that naturally come to mind:

\begin{enumerate}
\item[(SY)] If $\Omega$ has symmetries, does $D$ have the same symmetries?
(Note that if $\Omega$ and $D$ have a symmetry in common then $u$ will
also have this symmetry since it is uniquely determined by $\Omega$ and $D$.)
\item[(CX)] Assume $\Omega$ is convex. Is $D^c:=\Omega\setminus D$ convex?
 Is $D$ unique?
\item[(CN)] Is $D$ or $D^c$ connected? 
\item[(FB)] What is the regularity of the {\em free boundary} $\partial D$?
\end{enumerate}

We give partial answers to all of these questions. Some proofs, mainly
relating to (FB), and additional results can be found in the companion paper
[CGK]. 
Many open problems remain, see Section \ref{secopen}.

At this point, the reader is invited to look at Figures 1-3 for a first
impression.
\vspace{\baselineskip}

We now state our qualitative results.
As a general convention, constants only depend on the quantities
indicated as subscripts or in parentheses, unless otherwise specified.
Often we suppress the subscript $\Omega$.

First, as an easy consequence of Theorem \ref{thexist} one has:

\begin{theorem} \label{thtub}
Fix $\alpha>0$, $A>0$, and let $D$ be an optimal configuration.
\begin{enumerate}
\item[(a)] $D$ contains a tubular neighborhood of the boundary $\partial
\Omega$.
\item[(b)] If $\alpha<\abar_\Omega(A)$ then every connected component $D_0$ of 
the interior of $D$
hits the boundary, i.e. $\overline{D_0}\cap\partial\Omega \not=\emptyset$.
\end{enumerate}
In particular, if $\Omega$ is simply connected and $\alpha<\abar_\Omega(A)$
then $D$ is connected.
\end{theorem}
The number $\abar_\Omega(A)$ was defined above, see
\eqref{eqabardef0}.
The significance of the condition $\alpha<\abar_\Omega(A)$ is that it is
equivalent to $\Delta u < 0$ on $\Omega$.
 One always has 
$$\abar_\Omega(A) \geq \mu_\Omega.$$
Here and throughout the paper, $\mu_\Omega$ denotes the first
eigenvalue of the Dirichlet Laplacian on $\Omega$, and $\psi_\Omega$
the positive, $L^2$-normalized
eigenfunction:
$$ -\Delta \psi_\Omega = \mu_\Omega \psi_\Omega \quad\text{on }\Omega,
\quad\psi_{\Omega} = 0\quad\text{on }\partial \Omega.$$

Next, we consider the dependence of $\Lambda_\Omega$ and solutions $(u,D)$
on $\alpha$ and $A$.
Here it is convenient to formulate our problem also for $\alpha=0$,
as follows: If $\alpha=0$ then a solution
(unique in this case) is a pair $(\psi_\Omega,D)$
where $D$ is the sublevel set of $\psi_\Omega$ of area $A$.
(Since $\psi_\Omega$ is real analytic and non-constant, such $D$ exists
for every $A$ and is unique.)

 We will prove strict monotonicity and 
Lipschitz continuity of $\Lambda_\Omega$ in both parameters (Proposition
\ref{propparam}). Continuous dependence of optimal pairs $(u,D)$
on the parameters may be expected only at parameter values
where they are unique.
This is the case, in particular, if $\alpha=0$ 
or $A=0$ or $A=|\Omega|$;
in these cases $u=\psi_\Omega$, and the continuity is proved in [CGK].
 Here we only state the results. They are used only in the proof
 of Theorem \ref{thasmallcxfb}.
For example, we have the following:
\begin{theorem} \label{thDclose}
For $s\geq 0 $ let $[\Omega]^s = \{\psi_\Omega\leq s\}$, where
$\psi_\Omega$ is the positive $L^2$-normalized first eigenfunction 
of $-\Delta$ on $\Omega$.
Fix $A\in [0,|\Omega|]$ and choose $t_\Omega$ such that $|[\Omega]^{t_\Omega}|=A$.
Then for any $\delta>0$ there is $\alpha_0=\alpha_0(\delta,\Omega)$ such that
whenever $\alpha<\alpha_0$ and $D$ is an optimal configuration for $(\alpha,A)$
then 
$ |t-t_\Omega| < \delta $
and
$$ [\Omega]^{t_\Omega-\delta} \subset D \subset [\Omega]^{t_\Omega+\delta}.$$
\end{theorem}

We now address questions of {\em symmetry}.
First, we prove {\em symmetry preservation} in the presence of convexity:

\begin{theorem} \label{thsymmpres}
Assume that the domain $\Omega$ is symmetric and convex with respect
to the hyperplane $\{x_1=0\}$. In other words, for each $x'=(x_2,\ldots,x_n)$
the set 
\begin{equation} \label{eqsy0}
 \{x_1:\, (x_1,x') \in \Omega\}
\end{equation}
is either empty or an interval of the form $(-c,c)$.

Then for any solution $(u,D)$ both $u$ and $D$ are symmetric
with respect to  $\{x_1=0\}$,  $D^c$ is convex with respect
to $\{x_1=0\}$, and $u$ is decreasing in $x_1$ for $x_1\geq 0$.
\end{theorem}
For example, any solution in an elliptic region has a double reflection
symmetry, see Figure 1.
The principal tool here is Steiner symmetrization. See [K2] for an
overview on such methods. Theorem \ref{thsymmpres} easily
implies the following uniqueness result (the only case where we
can prove uniqueness!):

\begin{cor} \label{corball}
Let $\Omega = \{|x|<1\}$ be the ball.
Then there is a unique optimal configuration  $D$
for any $\alpha,A$, and $D$ is a shell region
$$ D= \{x:\, r(A) < |x| < 1\}. $$
\end{cor}

One of the most interesting phenomena studied in this paper is 
{\em symmetry breaking} for certain plane domains $\Omega$. That is,
an optimal configuration $D$ may have less symmetry than $\Omega$. We will
prove it for two types of domains: Thin annuli and dumbbells with narrow handle.
An annulus has rotational symmetry, a dumbbell has a reflection symmetry.

\begin{theorem} \label{thsymmann}
Fix $\alpha>0$ and $\delta\in (0,1)$. For $a>0$ let
$\Omega_a = \{x\in\R^2:\, a<|x|<a+1\}$. There exists $a_0=a_0(\alpha,\delta)$
such that whenever $a>a_0$ and $D$ is an optimal configuration for
$\Omega_a$ with parameters $\alpha$ and $A=\delta|\Omega_a|$ then
$D$ is not rotationally symmetric.
\end{theorem}
See Figure 2.
For dumbbells we prove a little more than symmetry breaking:

\begin{theorem} \label{thdumb} For $h\in (0,1)$ define the dumbbell
with handle width $2h$
$$ \Omega_h = B_1(-2,0) \cup ((-2,2)\times (-h,h)) \cup B_1 (2,0) $$
where $B_r(p) = \{x\in\R^2:\, |x-p|<r\}$.
Fix $\alpha>0$ and $A\in (0,2\pi)$. Then there is $h_0=h_0(\alpha,A)>0$ 
such that we have for $h<h_0$:
\begin{enumerate}
\item[(a)] Any optimal pair $(u,D)$ is not symmetric with respect
to the $x_2$-axis.
\item[(b)] If $A>\pi$ then for any optimal pair $(u,D)$ the complement
$D^c$ is contained in one of the lobes (i.e. one of the balls $B_1(\pm2,0)$).
\end{enumerate}
\end{theorem}
See Figure 3.
In fact, similar results hold for more general dumbbells.

As we remarked before, symmetry breaking implies non-uniqueness: For
example for a dumbbell the pair $(u',D')$ obtained from a solution 
$(u,D)$ by reflection in the $x_2$-axis will be a solution, and
different from $(u,D)$ by the theorem. 
\vspace{\baselineskip}

The following result on the regularity of the free boundary is proved
in [CGK]:

\begin{theorem} \label{thfb}
If $(u,D)$ is an optimal pair, $x\in\partial D$ and 
$\grad u(x) \not=0$ then $\partial D$ is a real analytic hypersurface near
$x$.
\end{theorem}
The difficulty is that $\chi_D$ is discontinuous at $x\in\partial D$, so $u$ is
not even $C^2$ there. That the level set $\{u=t\}$ has $C^\omega$ regularity
nevertheless is proved by introduction of suitable local coordinates (with
$u$ as one coordinate) and analysis of the resulting nonlinear elliptic
equation.

Similar arguments and continuity considerations for $\alpha$ near zero
allow us to give partial answers to problems (CX) and (FB):
\begin{theorem} \label{thasmallcxfb}
Suppose $\Omega$ is convex and has a $C^2$ boundary. Then there
is $\alpha_0(A,\Omega)>0$ such that for any $\alpha<\alpha_0$ and
any optimal configuration $D$, one has:
\begin{enumerate}
\item[(a)] $\partial D\cap\Omega$ is real analytic.
\item[(b)] $D^c$ is convex.
\end{enumerate}
\end{theorem}
\vspace{\baselineskip}

Problem (P) and generalizations of it (to higher eigenvalues
and to a maximization problem), but with fewer qualitative results,
were studied before in
[Kr], [CM], and [C] (where Theorem \ref{thsymmpres} is stated, but
the proof is incomplete since the case of equality in the
rearrangement inequalities is not addressed).

Problems similar to problem (M) (e.g. with $L^p$ potentials) were considered in
[AH], [Eg], [AHS], [CL], and [HKK]. 
\vspace{\baselineskip}

The paper is organized as follows:
In Section \ref{secbasics}, we prove Theorems \ref{thexist} and
\ref{thtub} and discuss the parameter dependence of $\Lambda_\Omega$.
Also, in Subsection \ref{subsecPM} we discuss the relation of 
problems (P) and (M).
In Section \ref{secsymm} we prove Theorems \ref{thsymmpres},
\ref{thsymmann}, and \ref{thdumb} on symmetry questions, and 
Corollary \ref{corball}.
In Section \ref{secfbcx} we prove Theorem \ref{thasmallcxfb}.
In Section \ref{secnum} we describe the numerical algorithm used.
In Section \ref{secopen} we state some open problems and conjectures.
Finally, we collect some standard facts about elliptic
PDEs in the Appendix.

\vspace{\baselineskip}

\begin{center} ACKNOWLEDGMENT \end{center}
\nopagebreak

We started this work while K.\ Kurata  was visiting the 
Erwin Schr\"{o}dinger International Institute for Mathematical Physics (ESI) in 
Vienna. K.\ Kurata is  partially supported by Tokyo Metropolitan University 
maintenance costs and by ESI and  wishes to thank Professor 
T.\ Hoffmann-Ostenhof for his invitation and  the members of 
ESI for their hospitality.
 D.\ Grieser was also at the ESI, and thanks L.\ Friedlander for inviting him. 
  M.\ Imai deeply appreciates helpful comments and heartful encouragement 
by Professor Teruo Ushijima in the University of Electro-Communications
(Tokyo). 
 We thank Professor 
M.\ Loss for his interest and some valuable comments. 
We thank Professor E.\ Harrell for his interest in our work and for 
informing us of the related works [AH], [AHS], and valuable discussions. 
S. Chanillo was supported by NSF grant DMS-9970359.
\section{Basic results} \label{secbasics}

\subsection{Existence and regularity. Proof
of Theorem \ref{thexist}}
We first prove {\em existence and regularity:}
The regularity statements in (a) hold for solutions of
equations
$$ -\Delta u + \rho u = 0$$
with $\rho$ bounded by standard elliptic theory, see for example
[GT, Theorem 8.29 and Corollary 8.36]. 

To prove existence, fix $\alpha$ and $A$, and write $\Lambda=\Lambda_\Omega(\alpha,A)$,
$\lambda(D) = \lambda_\Omega(\alpha,D)$ for simplicity.
Let $D_j$ be a minimizing sequence, i.e. $\lambda(D_j)\to \Lambda$
as $j\to\infty$.
Let $u_j\in H_0^1$ (all function spaces are defined on $\Omega$)
 be the positive $L^2$-normalized first
eigenfunction of $-\Delta+\alpha\chi_{D_j}$.
Since $\lambda(D_j)$ is bounded, the sequence $\{u_j\}$ is bounded
in $H^1_0$. Also, $\{\chi_{D_j}\}$ is bounded in $L^2$.
Therefore, we may choose a subsequence (again denoted $u_j,D_j$)
and $u\in H^1_0$, $\eta\in L^2$ such that
$u_j \weakto u$ in $H^1_0$ (weak convergence) and
$\chi_{D_j}\weakto\eta$ in $L^2$. This implies 
$u_j\to u$ (strongly) in $L^2$, $\chi_{D_j} u_j \weakto \eta u$
in $L^2$, and $\int_\Omega \eta = A$.
Now taking limits in the weak form of the eigenvalue equation
$$ \int_\Omega \grad u_j\cdot \grad\psi + \alpha \int_\Omega
    \chi_{D_j} u_j\psi = \lambda(D_j) \int_\Omega u_j\psi 
    \qquad \forall \psi\in H^1_0$$
we get
\begin{equation} \label{equeta}
-\Delta u + \alpha \eta u = \Lambda u \quad\text{ (weakly).}
\end{equation}
We have
$$ 0\leq\eta\leq 1 \quad\text{ a.e.}$$
since $0\leq\chi_{D_j} \leq 1$ for all $j$ and weak convergence preserves
pointwise inequalities a.e. (exercise!). Therefore, $u$ has the regularity
stated in (a).

It remains to prove that $\eta$ may be replaced by a characteristic function.
Since $\int_\Omega u^2=1$, \eqref{equeta} shows that
\begin{equation} \label{eqetaineq}
 \int_\Omega |\grad u|^2 + \alpha \int_\Omega \eta u^2 = \Lambda.
\end{equation}
Now the minimization problem
$$ \inf_{\substack{\eta:\int\eta=A \\ 0\leq\eta\leq 1}}
    \int_\Omega \eta u^2 $$
has a solution $\eta=\chi_D$ where $D$ is any set with $|D| = A$ and
\begin{equation} \label{eqDineq}
\{u<t\} \subset D \subset \{u\leq t\},
\quad t:=\sup\{s:|\{u<s\}|<A\}
\end{equation}
(compare the 'bathtub principle', Theorem 1.18 in [LL]).
Therefore, we get from \eqref{eqetaineq} 
$$\int_\Omega |\grad u|^2 + \alpha \int_\Omega \chi_D u^2 \leq \Lambda.$$
By definition of $\Lambda$ as a minimum, this must actually be an
equality, and $(u,D)$ is a solution.

(b) Let $(u,D)$ be any solution. Then it is obvious that
\eqref{eqDineq} must hold (always up to a set of measure zero;
if \eqref{eqDineq} didn't hold then one could reduce $\int_D u^2$ by
shifting a part of $D$ from $\{u>t\}$ to $\{u\leq t\}$).
Set $\calN_s = \{u=s\}$ for any $s>0$.
Using Lemma 7.7 from [GT] twice, we see that $\Delta u = 0$ a.e.
on $\calN_s$ (since $u\equiv$ const on $\calN_s$; recall that
$u$ is in $H^2$). Therefore,
\begin{equation} \label{equaeonN}
(\Lambda-\alpha\chi_D)u = 0 \quad\text{ a.e. on }\calN_s.
\end{equation}
Since $u>0$ and $\Lambda>0$, this shows that $D^c\cap\calN_s$ has
measure zero. Taking $s=t$ we get (b).

(c) If $s>t$ then $\calN_s\subset D^c$, so $|\calN_s|=0$ by 
\eqref{equaeonN}.
The same argument works if $s=t$ and $\alpha\not=\Lambda$.

Finally, $u$ satisfies $-\Delta u = (\Lambda-\alpha) u$
on the open set $\{u<t\}$, hence $u$ is real analytic there,
and therefore the level sets $\calN_s$ have measure zero for $s<t$.
\qed

\pf[Proof of Theorem \ref{thtub}]
Part (a) is clear from Theorem \ref{thexist}(b).
To prove (b), assume this was false. 
Then there is an open subset $D_0\subset \{u\leq t\}$ with
$\partial D_0 \subset \overline{D^c} = \{u\geq t\}$ and therefore
$u=t$ on $\partial D_0$.
Then $u$ assumes a minimum at some  $x_0\in D_0$.
But this is a contradiction since $\alpha<\abar_\Omega(A)$ implies
$\Lambda(\alpha,A) > \alpha$ (see Proposition \ref{propparam} below)
and therefore
$\Delta u = (\alpha-\Lambda(\alpha,A))u < 0$ on $D_0$.
\qed

\subsection{Parameter dependence of $\Lambda$} 

\begin{prop}  \label{propparam}
\begin{enumerate}
\item[(a)] The function $(\alpha,A)\mapsto \Lambda(\alpha,A)$
is Lipschitz continuous, uniformly on bounded sets. More precisely,
we have, for any $\alpha,\alpha'\geq 0$, $A,A'\in [0,|\Omega|]$,
\begin{multline} \label{eqlipest}
 |\Lambda(\alpha,A) - \Lambda(\alpha',A')|  \leq \\
  |\alpha-\alpha'|\,\frac{\max\{A,A'\}}{|\Omega|}
  + |A-A'|\, \min\{\alpha,\alpha'\} C_{\Omega,\max\{\alpha,\alpha'\}}
\end{multline}
with $C_{\Omega,\alpha}$ bounded for $\alpha$ bounded.

\item[(b)] $\Lambda(\alpha,A)$ is strictly increasing in $A$ for
fixed $\alpha>0$, strictly increasing in $\alpha$ for fixed $A>0$,
and $\Lambda(\alpha,A) - \alpha$ is strictly decreasing in $\alpha$
for fixed $A<|\Omega|$.
\item[(c)]
If $A<|\Omega|$ then there is a unique value $\alpha=\abar_\Omega(A)$ with
\begin{equation} \label{eqabardef}
\Lambda(\abar_\Omega(A),A) = \abar_\Omega(A).
\end{equation}
The function $\abar_\Omega$ is continuous and strictly increasing, 
$\abar_\Omega(0) = \mu_\Omega$ and $\abar_\Omega(A) \to\infty$
as $A\to |\Omega|$.
\end{enumerate}
\end{prop}

\pf
(a)
Write $\Lambda=\Lambda(\alpha,A)$ and $\Lambda'=\Lambda(\alpha',A')$,
and let $(u,D)$, $(u',D')$ be minimizers for $\Lambda$, $\Lambda'$
respectively. We may assume $\int_\Omega u^2 = \int_\Omega (u')^2 = 1$, so that
\begin{alignat*}{2}
\Lambda &= \int_\Omega |\grad u|^2 + \alpha \int_D u^2, & \quad |D| &= A,
\end{alignat*}
and similarly for $\Lambda'$ etc.
By symmetry of \eqref{eqlipest} we may assume
that $A'\geq A$. Choose $D_1\subset D'$ with $|D_1|=A$
and $D'_1\supset D$ with $|D'_1| = A'$. Here we may assume that $D'_1$
is of the form $\{u\leq s\}$ for a suitable number $s$.
Using the optimality of $(u,D)$ for $\Lambda$ we get
\begin{equation} \label{eqlamineq1}
\Lambda \leq \int_\Omega |\grad u'|^2 + \alpha \int_{D_1} (u')^2 
        = \Lambda' + (\alpha-\alpha')\int_{D'} (u')^2 - \alpha
                \int_{D'\setminus D_1} (u')^2.
\end{equation}
Similarly, using the optimality of $(u',D')$ for $\Lambda'$ we get
\begin{equation} \label{eqlamineq2}
\Lambda' \leq \int_\Omega |\grad u|^2 + \alpha' \int_{D'_1} u^2 
         = \Lambda + (\alpha'-\alpha)\int_{D'_1} u^2 + \alpha
                \int_{D'_1\setminus D} u^2.
\end{equation}
Alternatively, we may rewrite this as
\begin{equation} \tag{\ref{eqlamineq2}'}
\Lambda' \leq \Lambda + (\alpha'-\alpha)\int_D u^2 + \alpha'\int_{D'_1\setminus D}
                         u^2.
\end{equation}
In order to estimate the integrals in \eqref{eqlamineq1}, \eqref{eqlamineq2}
and (\ref{eqlamineq2}') which are multiplied by $\pm(\alpha-\alpha')$,
observe that for any $s>0$ and any function $u$ we have
$$ \frac {\int_{\{u\leq s\}} u^2}{\int_\Omega u^2} \leq \frac{|\{u\leq s\}|}{|\Omega|}.$$
The other integrals are estimated using the uniform estimate \eqref{eqpdeunif}:
$u$ solves the equation $-\Delta u + \alpha \chi_D u = \Lambda u.$
$\Lambda$ is bounded in terms of $\Omega$ and $\alpha$ since one
may apply \eqref{eqlamineq1} with $\alpha'=0, A=A'$, to obtain
$ \Lambda \leq \mu_\Omega + \alpha.$
Therefore, the uniform bound \eqref{eqpdeunif}, applied to $G=\Omega$, yields
$$ \int_{D'_1\setminus D} u^2 \leq (A'-A) \sup_\Omega u^2 \leq
(A'-A) C_{\Omega,\alpha}. $$
Finally, we obtain \eqref{eqlipest} by applying these estimates to
\eqref{eqlamineq1} and \eqref{eqlamineq2} in the case $\alpha\leq\alpha'$,
and to \eqref{eqlamineq1} and (\ref{eqlamineq2}') if $\alpha\geq\alpha'$.

(b) This follows immediately from \eqref{eqlamineq1} and the unique
continuation theorem.

(c) This follows easily from (a) and (b) since $\Lambda(\alpha,A)-\alpha$
equals $\mu_\Omega>0$ for $\alpha=0$ and tends to $-\infty$ as $\alpha\to
\infty$ by (a).
\qed

We now consider continuous dependence of optimal pairs $(u,D)$ on 
the data. First, near $\alpha=0$:

\begin{prop} \label{propparam2}
Fix $D\subset\Omega$. Let $u_{\alpha,D}$ be the (positive,
$L^2$-normalized) first eigenfunction of
$-\Delta + \alpha\chi_D$, and $\psi_\Omega=u_{0,D}$ the first 
eigenfunction of $-\Delta$.
Then there is a constant $C=C_\Omega$ such that, for $0\leq\alpha\leq 1$,
\begin{align*}
 \|u_{\alpha,D}-\psi_{\Omega}\| &\le C \alpha,
\end{align*}
in the $H^2(\Omega)$ and $L^\infty(\Omega)$ norms, and in $C^{1,\delta}(\Omega)$
if $\partial\Omega$ is in $C^{1,\delta}(\Omega)$.
\end{prop}

\pf 
See [CGK]. \qed

\pf[Proof of Theorem \ref{thDclose}]
This is almost immediate from Proposition \ref{propparam2}, see [CGK]. \qed

Similarly, one has continuity in $A$ at $A=0$ and at $A=|\Omega|$. Here we only
consider the latter case:
\begin{prop} \label{propAlarge}
Let $\Omega$ be a smooth bounded domain and fix $\alpha >0$. Let 
$$ M = \max_\Omega \psi_\Omega.$$
Then, for any $\delta >0$ there is $A_0=A_0(\delta,\alpha,\Omega)<|\Omega|$ 
such that whenever $A > A_0$ and $D$ is an
optimal configuration for $(\alpha,A)$ then
$$
D^c
 \subset \{\psi_\Omega > M -\delta\}.
$$ 
\end{prop}

\pf
See [CGK]. \qed

\subsection{Relation of problems (P) and (M)} \label{subsecPM}

We want to show that problem (P) 
(see Section \ref{secintro}) is a special case of problem (M).

The mathematical formulation of problem (P) is:

Given $0\leq h < H$ (lower and upper bounds for the densities of the materials
that 
are available) and the prescribed total mass $M\in [h|\Omega|,H|\Omega|], M>0$, 
consider measurable 'density functions'
$\rho$ satisfying
$$ h\leq \rho \leq H,\quad \int_\Omega \rho = M.$$
Then the objective is to find $\rho$ and $u$ which realize the minimum in
\begin{equation} \label{eqPvar}
\Theta (h,H,M) := \inf_\rho \inf_{u\in H_0^1(\Omega)}
  \frac {\int_\Omega |\grad u|^2} {\int_\Omega \rho u^2}.
\end{equation}
The corresponding eigenvalue problem is
\begin{equation} \label{eqPpde}
 -\Delta u = \Theta\rho u,\qquad u_{|\partial\Omega} = 0.
\end{equation}
(We assume the modulus of elasticity to be the same for all materials.)
Problem (P) and problem (M) are related in the following way:

\begin{theorem} \label{thprobPM}
\begin{enumerate}
\item[(a)] If $(u,\rho)$ is a minimizer for problem (P) then
 $\rho$ is of the form
 $$ \rho_D = h\chi_D + H\chi_{D^c}$$
 for a set $D$ of the form $D=\{u\leq t\}$. That is, only two types of
 materials occur.
\item[(b)] The pair $(u,\rho_D)$ is a minimizer for problem (P), 
 with parameter values $(h,H,M)$, if and only if $(u,D)$ is a 
 minimizer (optimal pair) for problem (M), with parameter values
 $(\alpha,A)$ given by
 \begin{align}
 \alpha &= (H-h) \Theta(h,H,M), \label{eqalphaP} \\
 A & = \frac{H|\Omega| - M} {H-h}.  \label{eqAP}
 \end{align}
 The minimal eigenvalues are related by
 \begin{equation} \label{eqLP}
  \Lambda(\alpha,A) = H\Theta(h,H,M).
 \end{equation}
\item[(c)] The values of $(\alpha,A)$ that occur when $h,H,M$ vary are
precisely those satisfying
\begin{alignat*}{2}
A &\in [0,|\Omega|), &\quad 0&<\alpha\leq \abar_\Omega(A) \quad \text{or} \\
A& = |\Omega|, &\quad 0&<\alpha<\infty,
\end{alignat*}
where $\abar_\Omega(A)$ is defined in \eqref{eqabardef}.
In particular, $\alpha=\abar_\Omega(A)$ corresponds to $h=0$.
\end{enumerate}
\end{theorem}
Note that problem (P) really depends on two parameters only since
for $\kappa>0$ one has $$\Theta(\kappa h,\kappa H,\kappa M) =
 \kappa^{-1} \Theta(h,H,M),$$
with the same minimizers (up to a factor $\kappa$ for $\rho$). 
This is obvious from \eqref{eqPvar}.

\pf
(a) This is almost obvious from \eqref{eqPvar}, and
 proved just like part (b) of Theorem \ref{thexist}.

(b) 
First, if $\rho=\rho_D$ and $|D|=A$ then
$ M = \int_\Omega \rho = Ah + (|\Omega|-A)H,$
which gives \eqref{eqAP}.

Simple manipulation shows that
\begin{equation}
 -\Delta u = \Theta\rho_D u = \Theta(h\chi_D + H\chi_{D^c}) u \label{eqPM1} 
\end{equation}
is equivalent to 
\begin{equation} 
     -\Delta u + (H-h) \Theta\chi_D u  = H\Theta u \label{eqPM2}.
\end{equation}
Now if $(u,\rho_D)$ is a minimizer for problem (P) then it satisfies
\eqref{eqPM1} with $\Theta = \Theta(h,H,M)$,
and then \eqref{eqPM2} shows that $\Lambda(\alpha,A)
\leq H\Theta(h,H,M)$ with $\alpha$ satisfying \eqref{eqalphaP}.

Conversely, if $(u,D)$ is a minimizer for problem (M) with parameter
values $(\alpha,A)$ given by \eqref{eqalphaP}, \eqref{eqAP}
then \eqref{eqPM2} holds with $H\Theta$ replaced by 
$\Lambda=\Lambda(\alpha,A)$, so instead of \eqref{eqPM1} we get
$-\Delta u = \Theta \rho_D u + (\Lambda-H\Theta)u$
where $\Theta=\Theta(h,H,M)$. Multiplying by $u$ and integrating gives
$$ \int_\Omega |\grad u|^2 = \Theta \int_\Omega \rho_D u^2 + (\Lambda-H\Theta)
  \int_\Omega u^2.$$
Now the definition of $\Theta$ implies that
$ \int_\Omega |\grad u|^2 \geq \Theta \int_\Omega \rho_D u^2,$
so we get $\Lambda\geq H\Theta$.
This proves $\Lambda(\alpha,A) = H\Theta(h,H,M)$ and part (b).

(c)
If $A=|\Omega|$ then $D=\Omega$, $\rho\equiv h$ and therefore
$h\Theta(h,H,M) = \mu_\Omega$ from \eqref{eqPpde},
so $\alpha = \frac{H-h}{h}\mu_\Omega$ can take any positive value
by suitable choice of $h$ and $H$.

Now let $A<|\Omega|$. By Proposition \ref{propparam}(b) and (c), $\alpha$ varies
in the indicated range precisely when $\Lambda(\alpha,A) - \alpha$
varies in $[0,\mu_\Omega)$. From \eqref{eqalphaP} and \eqref{eqLP} one has
$$\Lambda (\alpha,A) - \alpha = h\Theta:= h\Theta(h,H,M),$$
so we only need to show that $h\Theta$ has range  
$[0,\mu_\Omega)$ (with $A$ fixed).
First, $h\Theta\geq 0$ by definition, and $h\Theta = \Lambda-\alpha < \mu_\Omega$
by Proposition \ref{propparam}, since $\alpha = (H-h)\Theta > 0$, so the
range of $h\Theta$ is contained in $[0,\mu_\Omega)$.
Next, $h\Theta=0$ for $h=0$ (and then $M$ can be adjusted to $A$),
and in the limit $H=h$ one has $\rho\equiv h$ and $h\Theta = \mu_\Omega$,
so when $H\to h$ then $h\Theta\to\mu_\Omega$, and clearly $M$ can be adjusted
to $A$.
Using continuity of $h\Theta$ 
(which is proved as for $\Lambda$ in Proposition \ref{propparam})
 we get the claim.
\qed

\section{Symmetry preservation and symmetry breaking} \label{secsymm}

\subsection{Symmetry preservation in the presence of convexity}
Here we prove Theorem \ref{thsymmpres}.

\pf[Proof of Theorem \ref{thsymmpres}]
We use Steiner symmetrization (symmetrically decreasing rearrangement)
$u\mapsto u^\#$ 
with respect to the hyperplane $\{x_1=0\}$.
This is defined as follows. Assume $u\in H^1_0(\Omega)\cap C^0(\Omega)$: 
For each $x'$, $u^\#(\cdot,x')$ is the unique
function of $x_1$ which is symmetric in $x_1$ and decreasing for $x_1\geq 0$
such that $|\{x_1: u^\#(x_1,x') > t\}| = |\{x_1: u(x_1,x') > t\}|$
for all $t\in \R$.
It is well-known (see, e.g., [LL], [AB]) that, for all $x'$ and
$i=1,\ldots,n$, with integrals taken over the set \eqref{eqsy0},
\begin{eqnarray}
\int |\partial_{x_i} u^\#|^2\,dx_1 &\le & \int|\partial_{x_i} u|^2\,dx_1,  
  \label{eqsy1}\\
\int(u^\#)^2\,dx_1 &=& \int u^2\,dx_1, \label{eqsy1'} \\
\int(\alpha\chi_D)_\#(u^\#)^2\,dx _1
   &\le& \int\alpha\chi_D u^2\,dx_1. \label{eqsy2}
\end{eqnarray}
Here,  $f_\#$ is the increasing symmetric rearrangement of a function $f$, 
which is defined by $f_\#=-(-f)^\#$.
Note that \eqref{eqsy1} for $i=1$ is just the standard rearrangement
inequality in one dimension, while for $i>1$ it is proved as follows:
Replace the partial derivatives by difference quotients
$(v_\eps(x_1) - v_0(x_1))/\eps$
with $v_\eps(x_1) = u(x_1,\ldots,x_i+\eps,\ldots)$.
After multplication by $\eps^2$ the claimed inequality 
becomes simply $ \int |v_\eps^\#-v_0^\#|^2 dx_1 \leq \int |v_\eps-v_0|^2 dx_1$
which is well-known.

Fix $\alpha$ and $A$ and assume $(u,D)$ is an optimal pair.
Define the set $D^\#$ by $\chi_{D^\#} = (\chi_D)_\#$. 
Integrating \eqref{eqsy1}, \eqref{eqsy1'} and \eqref{eqsy2}
over $x'$ and summing \eqref{eqsy1} over $i$ we get
\begin{eqnarray}
\lambda(\alpha,D^\#) &\le&\frac{\int_{\Omega}|\nabla u^\#|^2\,dx + 
 \int_{\Omega} (\alpha \chi_D)_\#(u^\#)^2\,dx}
{\int_{\Omega}(u^\#)^2\,dx}\nonumber\\
&\le&
\frac{\int_{\Omega}|\nabla u|^2\,dx 
+ \int_{\Omega} \alpha \chi_D u^2\,dx
}{\int_{\Omega}u^2\,dx}
= \lambda(\alpha;D).
\end{eqnarray}
Since we have $|D^\#| = |D|=A$ (by \eqref{eqsy1'} applied to $\chi_D$),
optimality of $(u,D)$ implies that $(u^\#,D^\#)$ is also a minimizer
and that equality holds in \eqref{eqsy1} and \eqref{eqsy2}, for all $i$
and almost all $x'$.
We need to show that this implies $u=u^\#$. The statements about $D$
then follow from the characterization $D=\{u\leq t\}$.

First note that since $(u^\#,D^\#)$ is a minimizer, the function $u^\#$ solves
the equation
$-\Delta u^\# +\alpha\chi_{D^\#} u^\#=\lambda(\alpha;D^\#)u^\#$.
Therefore, $u$ and $u^\#$ are continuously differentiable by Theorem
\ref{thexist}, so equality in  \eqref{eqsy1} holds for all $x'$.
By a result of Brothers and Ziemer (see [BZ]) 
this equality implies $u^\#(x_1,x')= u(x_1,x')$ 
 for all $x_1$ provided the set
$ \{x_1: \partial_{x_1} u^\#(x_1,x') = 0\}$ has measure zero.

Therefore, we will be done once we have shown that 
the set 
\begin{equation}
\{v = 0\}\quad\text{ has measure zero, where }
v=\partial_{x_1} u^\#. \tag{*}
\end{equation}
We will give two proofs of this: The first proof works whenever
$\alpha\not=\abar_\Omega(A)$ and the second proof works whenever
$\alpha\leq\abar_\Omega(A)$, so together they cover all cases.

First proof of (*), assuming $\alpha\not=\abar_\Omega(A)$:
 Assume this was not so. Define $t^\#$ by $D^\#=\{u^\#\leq t^\#\}$. 
$v$ satisfies $-\Delta v + \alpha\chi_{D^\#} v= \lambda(\alpha,D^\#) v$
on $\{u^\#\not=t^\#\}$.
Since $\{u^\#=t^\#\}$ has measure zero by Theorem
\ref{thexist} and the assumption $\alpha\not=\abar_\Omega(A)$, 
$v$ vanishes on a set of positive
measure in the open set $\{u^\#\not= t^\#\}$, 
so the unique continuation theorem
(for sets of positive measure, see [FG])
applied to $v$ implies that $v\equiv 0$ on some connected component $K$
of $\{u^\#\not= t^\#\}$. Therefore, $u^\#$ is constant in
the $x_1$-direction on $K$. 
Since $u^\#=0$ or $t^\#$ on $\partial K$ 
we conclude that then $u^\#$ must actually
be constant on $K$. This is a contradiction to Theorem \ref{thexist}(c).

Second proof of (*), assuming $\alpha\leq\abar_\Omega(A)$
(this proof is taken from Cox [C]): We show that actually $v<0$
for $x_1>0$, so that $\{v=0\}$ is contained in the hyperplane
$\{x_1=0\}$.
We have $-\Delta u^\# = \Lambda(\alpha,A) u^\# - \alpha \chi_{D^\#} u^\#$, and
the right hand side is decreasing in $x_1$ (for $x_1> 0$) by definition
of the rearrangement and since $\alpha\leq\Lambda(\alpha,A)$ by
Proposition \ref{propparam}. Taking the $x_1$-derivative (in the sense
of distributions), we get $\Delta v \geq 0$ as distribution. Also,
$v$ is continuous, so by the classical theory of subharmonic functions
it satisfies the maximum principle (alternatively, it is in $H^1$ and
then the maximum principle as in [GT], Ch. 8, applies). Since
$v\leq 0$, we conclude that $v<0$ unless $v$ vanishes identically
in $x_1>0$, which is clearly impossible. This proves (*).

This concludes the proof that $u=u^\#$ and hence the proof of the theorem.
Note that in the case $\alpha\leq\abar_\Omega(A)$ the second proof of (*)
above actually shows that $u_{x_1}<0$ for $x_1>0$.

\qed

\pf[Proof of Corollary \ref{corball}]
The only set $D\subset \{|x|<1\}$ which has the symmetry and convexity
properties stated in Theorem \ref{thsymmpres} in all directions
is a shell region as stated. Clearly, $r(A)$ is uniquely determined
by $A$. Therefore, $D$ is unique. 
\qed

\subsection{Symmetry breaking on annuli}

We now give the proof of Theorem \ref{thsymmann} about 
symmetry breaking on an annulus, 
$$
\Omega=\Omega_a=\{x\in {\bf R}^2; a < |x| < a +1\},\quad a> 0.
$$
Let  $D$ be any radial set in $\Omega$, 
$$D=\{ (r,\theta ) ; r\in D_1, 0\le \theta < 2\pi\},\quad D_1\subset (a,a+1),$$
and let $u$ be the first eigenfunction for $D$, with eigenvalue $\sigma$:
\begin{equation} \label{eqannpde}
-\Delta u +\alpha \chi_D u=\sigma u \quad \text{ on } \Omega, 
\quad u|_{\partial \Omega}=0.
\end{equation}
For $a$ sufficiently large (depending on $\alpha$ and $\delta=|D|/|\Omega|$)
we 
will construct a comparison domain $\Dtilde$ and a function $\utilde$
which satisfy
\begin{equation} \label{eqanngoal}
\frac{\int_{\Omega_a} |\grad \utilde|^2 + \int_{\Omega_a} \chi_{\Dtilde} \utilde^2}
{\int_{\Omega_a} \utilde^2} \overset{!}{<} \sigma.
\end{equation}
This shows that $D$ is not an optimal configuration and hence
implies the theorem.

In order to construct $\Dtilde$ and $\utilde$, first pick $N=N(\delta)$
with
$$\delta < 1-\frac1{2N}$$
and consider the sector
$$
E_+= \Omega_a\cap \{ (r,\theta);  0\le \theta \le \pi/N\}.
$$
Then let $\utilde$ be the first Dirichlet eigenfunction of the Laplacian
on $E_+$ and $\lambda_1(E_+)$ be the first eigenvalue, 
\begin{equation} \label{eqannpde2}
-\Delta \utilde = \lambda_1(E_+) \utilde \quad\text{ on } E_+,
\quad \utilde|_{\partial E_+} = 0,
\end{equation}
 extended by zero on $\Omega\setminus E_+$; the set $\Dtilde$ can
be taken to be any subset of $\Omega\setminus E^+$ with $|\Dtilde|=|D|$.
This is possible since 
$|D|/|\Omega| = \delta < 1-\frac1{2N} = |\Omega\setminus E_+|/|\Omega|$.

Note that since $\supp \utilde \cap \Dtilde=\emptyset$, we have
$(\int_{\Omega_a} |\grad \utilde|^2 + \int_{\Omega_a} \chi_{\Dtilde} \utilde^2)/
\int_{\Omega_a} \utilde^2 = \int_{E_+} |\nabla \utilde|^2/
\int_{E_+} \utilde^2 = \lambda_1(E_+)$, so \eqref{eqanngoal} is equivalent to
\begin{equation} \label{eqanngoal1}
\lambda_1(E_+) \overset{!}{<} \sigma.
\end{equation}

In order to prove this,
we need to introduce a third eigenvalue problem, 
which is intermediate
between \eqref{eqannpde} and \eqref{eqannpde2}.

Define $v$ to be the lowest eigenfunction for the problem \eqref{eqannpde}
among functions of the form
$$ v(r,\theta) = h(r)\sin N\theta,$$
and let $\tau$ be the associated eigenvalue.
Note that problem \eqref{eqannpde} for such functions is equivalent to the
problem
\begin{gather} \label{eqhode}
-h''(r)-\frac{1}{r}h'(r)+\frac{N^2}{r^2}h(r)
+\alpha\chi_{D_1}(r)h(r)=\tau h(r) \quad \text{ on } [a,a+1],\\
h(a)= h(a+1)=0
\end{gather}
for $h$. Thus, $h$ is the first eigenfunction of this Sturm-Liouville problem,
and the eigenvalue $\tau$ is characterized by
\begin{equation} \label{eqhvar}
\tau =\inf_{g\in {\calS}}\frac{\int_a^{a+1}
((g')^2 + (\alpha\chi_{D_1} + \frac{N^2}{r^2})g^2 )r\,dr }
{\int_a^{a+1}g^2r\,dr},
\end{equation}
where ${\calS}=\{ g\in C^1[a, a+1]; g(a)=g(a+1)=0\}$. From this
the (well-known) fact that $h$ does not change sign on $[a,a+1]$
is evident; so we may assume $$h\geq 0.$$
We will compare $u$ with $v$ and $v$ with $\utilde$. The following
two lemmas provide the needed estimates.

\begin{lem} \label{propsymma}
Let $\sigma$ be the lowest eigenvalue for the problem \eqref{eqannpde}
(with $D$ radial)
on $\Omega_{a,b} = \{x\in \R^2: a<|x|<b\}$, and let $\tau$ be the lowest 
eigenvalue for eigenfunctions of the form 
$v(r,\theta)=h(r)\sin N\theta$ on $\Omega_{a,b}$. Then  we have
$$
 \tau-\sigma \le N^2/a^2.
$$
\end{lem}

\pf 
Since $\chi_D$ is assumed radial, the first eigenfunction of \eqref{eqannpde}
is a radial function $u=f(r)$. Now consider the trial function
$w(r,\theta)=f(r)\sin N\theta$. 
We have
$$
\tau \le \frac{\int_{\Omega_{a,b}}(|\nabla w|^2 + \alpha\chi_D w^2)\,dx}
{\int_{\Omega_{a,b}} w^2\,dx}.
$$
Thus, 
$$
\tau\le 
\frac{\int_a^b( (f'(r))^2 + \frac{N^2}{r^2}f(r)^2 
+ \alpha \chi_{D_1} f(r)^2)r\,dr }
{\int_a^bf(r)^2r\,dr}.
$$
By definition of $f(r)$ we get
$$
\tau \le \sigma +
\frac{\int_a^b(\frac{N^2}{r^2}f(r)^2)r\,dr}
{\int_a^bf(r)^2r\,dr}\le \sigma + N^2/ a^2.
$$
The claim follows.
\qed

\begin{lem} \label{lemsymma}
Define $v$ as above. Assume $D$ is radial and
$|D|/|\Omega|=\delta$. There exists a positive constant 
$c_{\alpha,\delta}$, independent of $a$, 
such that for all $a \ge 1$ we have
$$
\frac{\int_Dv^2\,dx}{\int_{\Omega}v^2\,dx}\ge c_{\alpha,\delta}.
$$   
\end{lem}

\pf 
We see from $v(r,\theta)=h(r)\sin N\theta$ that
\begin{equation}
\frac{\int_Dv^2\,dx}{\int_{\Omega}v^2\,dx}
=
\frac{\int_a^{a+1}\chi_{D_1}(r)h(r)^2r\,dr}
{\int_a^{a+1}h(r)^2r\,dr}.
\label{eq:9.2}
\end{equation}
$h$ satisfies equation \eqref{eqhode}. For $\tau$ one
has a uniform bound $\tau\leq C_{\alpha,\delta}$ with 
$C_{\alpha,\delta}$ independent of $a\geq 1$, because
from \eqref{eqhvar} one gets 
$$ \tau \leq \inf_{g\in {\calS}} \frac {\int_a^{a+1} (g')^2 r\,dr}
                                         {\int_a^{a+1} g^2 r\, dr}
                          + \alpha + N^2,$$
and by using for $g$ the translate of any fixed test function on $[0,1]$ 
one sees that the first term on the right is bounded by some absolute
constant.

Therefore, the coefficients of equation \eqref{eqhode} are 
uniformly bounded for $a\geq 1$. Also, we have $h\geq 0$.
Lemma \ref{lhar} in Section \ref{secpdefacts} then implies that
one has
\begin{equation} \label{eqh1}
\inf_{[a+\delta/4,a+1-\delta/4]} h \geq c_{\alpha,\delta} \|h\|_{L^2(a,a+1)}.
\end{equation}
Since $|D_1|=\delta$, we
have $|[a+\delta/4,a+1-\delta/4] \cap D_1| \geq \delta/2$. Therefore,
\begin{equation} \label{eqh2}
\int_a^{a+1} \chi_{D_1}(r)h(r)^2r\, dr \geq 
          \frac\delta2 a \inf_{[a+\delta/4,a+1-\delta/4]} h^2
\end{equation}
and
\begin{equation} \label{eqh3}
\int_a^{a+1} h(r)^2r\, dr \leq (a+1) \int_a^{a+1} h^2 \leq 2a\int_a^{a+1}h^2.
\end{equation}
Combining \eqref{eqh1}, \eqref{eqh2} and \eqref{eqh3} with 
\eqref{eq:9.2} we get the Lemma.
\qed

\pf[End of proof of Theorem \ref{thsymmann}]
We have
\begin{equation} \label{eqanntau}
\tau= 
 \frac{\int_{\Omega}|\nabla v|^2\,dx}{\int_{\Omega}v^2\,dx}
+\frac{\alpha\int_{\Omega}\chi_Dv^2\,dx}{\int_{\Omega}v^2\,dx}.
\end{equation}
Since $v(r,\theta)=h(r)\sin N\theta$, $v$ vanishes on the rays 
$\theta=0$ and $\theta=\pi/N$. Since $|v|$ and $|\grad v|$ are periodic
in $\theta$ of period $\pi/N$, we can replace $\Omega$ by $E_+$
in the first quotient.
 Therefore, we can use $v$ as test function in the Rayleigh quotient for
 the Dirichlet Laplacian on $E_+$ and obtain
$$ \frac{\int_\Omega |\grad v|^2 \, dx} {\int_\Omega v^2 \, dx}
 = \frac{\int_{E_+} |\grad v|^2 \, dx} {\int_{E_+} v^2 \, dx}
 \geq \lambda_1(E_+).$$
Combining this with \eqref{eqanntau} and Lemma \ref{lemsymma} we therefore get
\begin{equation} \label{eqsya1}
\tau\ge  \lambda_1(E_+) + \alpha c_{\alpha,\delta}.
\end{equation}
From Lemma \ref{propsymma} we then get
$$ \sigma > \tau - N^2/a^2 \geq \lambda_1(E_+) + 
\alpha c_{\alpha,\delta}-N^2/a^2.$$
If $a$ is chosen so large that $N^2/a^2 \leq \alpha c_{\alpha,\delta}$
then this gives \eqref{eqanngoal1} and hence the theorem.
\qed

\subsection{Symmetry breaking on dumbbells}

\pf[Proof of Theorem \ref{thdumb}]
Since $\alpha$ is fixed throughout, we will write $\lambda_\Omega(D)=
\lambda_\Omega(\alpha,D)$, $\Lambda_\Omega(A) = \Lambda_\Omega(\alpha,A)
= \inf_{|D|=A} \lambda_\Omega(D)$. Here we keep the index $\Omega$
since we will also consider these quantities with $\Omega$ replaced
by one of the \lq lobes' $B_\pm = B_1(\pm 2,0)$. All (implied) constants will
only depend on $\alpha$ and $A$. Write $\Lambda_B = \Lambda_{B_\pm}$,
and given $D$, let
$$ D_\pm = D\cap B_\pm,\quad A_\pm = |D_\pm|.$$
\newcommand {\Amin} {A_{\rm min}}
Further, we introduce
$$ \Amin  = \min\{\min(|D_-|,|D_+|):\, D\subset \Omega,
|D|=A\}.$$
Thus, if $D$ is distributed over $\Omega$ with the greatest possible
imbalance between $D_+$ and $D_-$
then the smaller of $D_\pm$ will have area $\Amin$.
It is easily checked that $$\Amin = \max (0,A-|B_-^c|).$$

We first sketch the idea of the proof:

1. For $h=0$, i.e. two disconnected balls, one clearly has
\begin{equation} \label{eqdumb0}
\Lambda_\Omega (A)= \min (\Lambda_B(A_-), \Lambda_B(A_+)).
\end{equation}
Since $\Lambda_B$ is strictly increasing,
it is optimal to put
as much of $D$ as possible in one ball, say $B_+$, and the \lq small' remainder
in the other. Thus $$\Lambda_\Omega(A) = \Lambda_B(\Amin),$$
and the eigenfunction is zero in $B_+$.

2. For small positive $h$, this situation should be  approximately
the same: Equation \eqref{eqdumb0} will hold with an error that is
a power of $h$ (compare equation (\ref{eqdumb2A}) below), so the
same argument as in  1. implies symmetry breaking. Also, the eigenfunction
must be small on one lobe, and since $D=\{u\leq t\}$, one gets
(b) from an estimate of $t$.
\vspace{\baselineskip}

We now carry out the details. 
Let $(u,D)$ be an optimal pair. Assume $$\int_\Omega u^2=1.$$
First we need an estimate ensuring that
the perturbation introduced by the handle is small.
This is provided by the following estimate near the boundary
(see [GT, Theorem 8.27 with $R_0=1$ and $R=3h$]), which is applicable
since $\Omega$ satisfies a uniform exterior cone condition (uniformly in $h$):
There is $\beta\in (0,1]$ such that
\begin{equation} \label{eqbdest}
\max_{x:\dist(x,\partial\Omega)\leq 3h} u(x) \leq Ch^\beta \|u\|_{L^2(\Omega)}.
\end{equation}

From this it follows that there is 
a cut-off function $\sigma = \sigma_h$ on $\Omega$ having the following
properties:
\begin{enumerate}
\item $0\leq\sigma\leq 1$ on $\Omega$. 
\item $\supp \sigma \subset B_-\cup B_+$.
\item $|u| = O(h^\beta)$ on $\supp (1-\sigma)$.
\item $\int_\Omega |\grad\sigma|^2 < C$, uniformly as $h\to 0$.
\end{enumerate}
To construct $\sigma$, choose $\chi\in C_0^\infty([0,2))$, $0\leq\chi\leq1$,
that equals one on $[0,3/2]$ and set $\sigma(x) = 1-\chi(|x-(\pm1,0)|/h)$ on
$B_\pm$ and $\sigma\equiv 0$ on the handle.
Properties 1,2 and 4 are easily checked directly, and property 3 follows
from \eqref{eqbdest}.

For brevity, denote, for $\Omega'\subset\Omega$,
$$Q_{\Omega'}(u) = \|\grad u\|_{L^2(\Omega')}^2 + 
           \alpha \|u\|_{L^2(D\cap\Omega')}^2$$
so that $ \Lambda_\Omega(A) = Q_\Omega(u).$
Without loss of generality we may assume
$$\Lambda_B(A_-) \leq \Lambda_B(A_+).$$

First, we show
\begin{equation} \label{eqdumb1}
\Lambda_B(\Amin) \geq \Lambda_\Omega(A).
\end{equation}
This is easy: Take an optimal pair $(\utilde,\Dtilde)$
for $\Lambda_B(\Amin)$,
extend $\utilde$ by zero outside $B_-$ and 
define a domain $\Dbar = \Dtilde\cup D'\subset\Omega$, $|\Dbar|=A$, by 
choosing any $D'\subset B_-^c$ with $|D'| = A-\Amin$.
Since $\utilde\equiv0$ on $D'$, one gets \eqref{eqdumb1} by
using $(\utilde,\Dbar)$ as 
a test pair for $\Lambda_\Omega(A)$.

Next, we show a reverse inequality. 
Using the properties of $\sigma$ and 
$\supp \grad\sigma \subset \supp (1-\sigma)$ 
we obtain, with $\|\cdot\|$ denoting the $L^2$-norm on $B_\pm$,
$\|\grad(\sigma u)\|^2 = \|\sigma\grad u + (\grad\sigma) u\|^2
\leq (\|\grad u\| +\|\grad\sigma\| \max_{\supp (1-\sigma)} u)^2 \leq 
\|\grad u\|^2 + O(h^\beta)$ and therefore
$Q_{B_\pm} (\sigma u) \leq Q_{B_\pm} (u) + O(h^\beta).$
Now we can use $\sigma u$ as test function for the lowest eigenvalue
of $-\Delta + \alpha \chi_{D\cap B_\pm}$ on $B_{\pm}$, and this gives the
third inequality in
\begin{eqnarray}
\Lambda_\Omega(A)=Q_\Omega(u) &\geq& \sum_\pm Q_{B_\pm}(u) \notag
            \geq \sum_\pm Q_{B_\pm}(\sigma u) - O(h^\beta)\notag\\
            &\geq& \sum_\pm \lambda_{B_\pm}(D_\pm) \int_{B_\pm} (\sigma u)^2 
                     - O(h^\beta)\notag\\
            &\geq&  \sum_\pm \lambda_{B_\pm}(D_\pm) \int_{B_\pm} u^2 
                     - O(h^\beta)      \notag \\
            &\geq& \Lambda_B(A_-) + (\Lambda_B(A_+)-\Lambda_B(A_-))
                         \int_{B_+}u^2  - O(h^\beta).  \label{eqdumb1a}
\end{eqnarray}                     
In the last two inequalities we have used 
property 3. of $\sigma$, the optimality of $\Lambda_B(A_\pm)$,
and $\int_\Omega u^2 = 1$.

Since we assume $\Lambda_B(A_+)\geq\Lambda_B(A_-)$, 
this and inequality (\ref{eqdumb1}) imply
\begin{equation} \label{eqdumb2A} 
\Lambda_B(\Amin) \geq \Lambda_\Omega(A) \geq \Lambda_B(A_-) - O(h^\beta).
\end{equation}
By strict monotonicity of $\Lambda_B$ one easily gets from this
$A_- \leq \Amin + o(1) \quad (h\to0).$

Next, from $D\subset D_+\cup D_- \cup H$ and $|H|<4h$ we have
$A< A_+ + A_- + 4h$, so
$A_+ - A_- > A-2A_- - 4h \geq A-2\Amin - o(1)$, and then
$\Amin = \max(0,A-|B_-^c|) \leq \max(0,A-\pi)$ gives
\begin{equation}\label{eqdumbb}
A_+ - A_- \geq \min (A,2\pi-A) - o(1).
\end{equation}
This shows $A_+\not= A_-$ for $h<h_0(A,\alpha)$ and therefore
proves part (a) the theorem.

Now we prove part (b). From \eqref{eqdumbb} we have $A_+ - A_- > c_0$
for some constant $c_0>0$, whenever $h<h_0(A,\alpha)$, so strict monotonicity and
continuity of $\Lambda_B$ imply
\begin{equation} \label{eqdumb2c}
 \Lambda_B(A_+) - \Lambda_B(A_-) > c
\end{equation}
with $c>0$ independent of $h$.
Now from \eqref{eqdumb1} and \eqref{eqdumb1a}, and using
$\Lambda_B(A_-)\geq\Lambda_B(\Amin)$ (since $A_-\geq \Amin$)
and monotonicity, we conclude
$(\Lambda_B(A_+) - \Lambda_B(A_-)) \int_{B_+} u^2  = O(h^\beta).$
This and \eqref{eqdumb2c} give
 $\int_{B_+} u^2 = O(h^\beta).
$
Since, by \eqref{eqbdest}, 
$u_{|\partial B_+} = O(h^\beta)$, this $L^2$ bound
implies a pointwise bound for $u$ on $B_+$ by \eqref{eqpdeunif}.
Combined with \eqref{eqbdest}, applied on the handle, this gives 
\begin{equation}\label{eqdumb4}
\sup_{x\not\in B_-} u(x) = O(h^{\beta/2}).
\end{equation}

Finally, we want to deduce from \eqref{eqdumb4} that $D^c\subset B_-$
if $A>\pi$ and $h$ is sufficiently small:
Since $(u,D)$ is an optimal pair, we have $D=\{u \leq t\}$ for
some $t>0$. Equation (\ref{eqdumb4}) shows that we are done if we can show 
that $t>c$ for
a constant $c>0$ independent of $h$. 

For $r\in (0,1)$  let $B_-(r)$ be the closed ball
of radius $r$ concentric with $B_-$. 
Applying Lemma \ref{lhar} to $G=B_-$ we see, since $\|u\|_{L^2(B_-)}
\geq 1 - O(h^\beta)$ by \eqref{eqdumb4}, that
\begin{equation}\label{eq6} \inf_{B_-(r)} u \geq c_r
\end{equation}
for any $r\in (0,1)$, with $c_r>0$ only depending on $r$, $A$ and $\alpha$,
 and this implies
$ |\{u \geq c_r\}| \geq |B_-(r)|.$
Therefore, we can conclude $t>c_r$ as soon as 
$|B_-(r)| > |\Omega| - A$. Since $|\Omega| \leq 2\pi + 4h$
and $A>\pi$, one can find such an $r$
if $h<h_0$, both $r$ and $h_0$ only
depending on $A$ (and $\alpha$). This completes the proof of the theorem.
\qed

\section{Free boundary and convex domains} \label{secfbcx}
\pf [Proof of Theorem \ref{thasmallcxfb}, Part (a)]
First recall, as a consequence of  results by Brascamp-Lieb [BL]  and 
Caffarelli-Spruck [CS], that the first eigenfunction $\psi$
on a convex domain possesses only
one point where $\nabla \psi=0$. This point is necessarily the point where
$\psi$ attains its maximum.
Now given $A$, we select $t_\Omega$ as in Theorem \ref{thDclose}, and we select 
$\delta_0<t_\Omega$
such that $t_\Omega+\delta_0< M$ where $M=\max_\Omega \psi$. With 
this choice of $\delta_0$ we use Theorem \ref{thDclose}
to determine a value $\alpha_1$ 
for which $[\Omega]^{t_\Omega-\delta_0} \subset 
 D\subset [\Omega]^{t_\Omega+\delta_0}$ for all $\alpha<\alpha_1$. Then
 the free boundary $\{u=t\}$ is contained in the closed annulus
 $A = \{t_\Omega - \delta_0 \leq \psi \leq t_\Omega + \delta_0\}$.
 We have $\grad\psi\not=0$ on $A$, so $C:=\min_A |\grad\psi|$ is positive.
Thus decreasing $\alpha_1$ 
to a smaller value $\alpha_0>0$, we can use Proposition \ref{propparam2} 
to conclude
that for all $\alpha<\alpha_0$ we have $|\nabla u| > C/2$
on $A$ and hence on the free boundary
$\{u=t\}$.
Applying Theorem \ref{thfb} we now get the first part of Theorem 
\ref{thasmallcxfb}.

\pf [Proof of Theorem \ref{thasmallcxfb}, part (b)]
We only sketch the proof.
Fix $x_0$ with $\grad\psi(x_0)\not=0$. 
Choose coordinates in which $\grad\psi(x_0)=(0,\ldots,0,a), a>0$, and 
for $x'$ near $x_0'$ (where $x'=(x_1,\ldots,x_{n-1})$) and $t$ near $t_0=\psi(x_0)$
denote
the locally unique solution $x_n$ of the equation $\psi(x',x_n)=t$
by $F_0(x',t)$. For $\alpha$ near zero and $x$ near $x_0$ one has
$\partial u_\alpha/\partial x_n\not=0$ by Proposition \ref{propparam2},
so we may define $F_\alpha$ similarly for $u_\alpha$ instead of $\psi$.

By a result of
Korevaar and Lewis [KL] the level set of $\psi$ through $x_0$ is
strictly convex, in the sense  that the matrix
$(\frac{\partial^2 F_0}{\partial x_i\partial x_j})_{i,j=1,\ldots,n-1}$
is positive definite at $(x_0',t_0)$. Therefore, the result
follows if one can show continuity of 
$\frac{\partial^2 F_\alpha}{\partial x_i\partial x_j}$ in $\alpha$
and $(x',t)$. Now the equation for $u$ gives for $F_\alpha$ a uniformly elliptic,
quasi-linear equation (writing $y=(x',t)$)
$$\sum_{i,j=1}^n b_{ij}(\grad F_\alpha)
     \frac{\partial^2 F_\alpha(y)} {\partial y_i \partial y_j} =
     \alpha\chi_{G_\alpha}(y_n)y_n-\Lambda(\alpha,A)y_n\ \ \  $$
with $b_{ij}$ real analytic
and  $G_\alpha=(-\infty,t_\alpha]$, where $t_\alpha$ is such that
$|\{u_\alpha\leq t_\alpha\}|=A$.
From this it is easy to derive the desired regularity, cf. the proof
of Lemma 3 in [CGK].
\qed 
\section{Numerical results} \label{secnum}   

     In this section we
     make a few remarks on our method for the numerical 
   solution of our eigenvalue problem.
   
   We use the finite element method for the discretization of 
   our eigenvalue problem, with conforming P-1 elements.
   To create the mesh we have utilized the 
   automatically spatial meshing program 
   encoded by Y. Tsukuda (see [TK]).
    In order to calculate the approximate first eigenvalue and the 
   corresponding eigenfunction, 
   we employ the power method. 

    Our method to obtain an optimal configuration
    is based on an algorithm that was introduced 
   in [Pi]. However, we do not insist on $D$ (the sought-for
   optimal configuration) to be a union of elements. This flexibility
   allows us to find a good approximation even without remeshing.
  
   We now describe the main procedure. 
   The given data are $A$ and $\alpha$.
   We first take any initial domain $D_0$ satisfying $|D_0| = A$. 
   Next, if we have obtained $D_{n-1}$ ($n=1,2,3,  \cdots$) then 
   we calculate the first eigenvalue 
   $\lambda_{n-1}$ and the corresponding eigenfunction $u_{n-1}$ 
   for the finite element approximation problem for the operator
   $-\Delta + \alpha\chi_{D_{n-1}}$. 
   Then we obtain $D_n$ from $u_{n-1}$ by finding a number $t_0$ such
   that $|\{u_{n-1}\leq t_0\}| = A$ and setting 
   $$D_n = \{u_{n-1} \leq t_0\}.$$
   The number $t_0$ is determined by a bisection method, i.e. by 
   setting $down_0=0$, $up_0=\max_\Omega u_{n-1}$, $j=0$ and
   then iterating Steps 1 and 2 (with
       $L(t) := |\{u_{n-1} \leq t \} |$)
       \begin{enumerate}
\item[Step 1:] Let $interm_j$: = $(up_j+down_j)/2$ and calculate $L(interm_j)$.
\item[Step 2:] If $L(interm_j)$ $<$ $A$, then $up_{j+1} := up_j$ and 
	 $down_{j+1} := interm_j$, else if $L(interm_j)$ $>$ $A$, 
	 then $up_{j+1} := interm_j$ and $down_{j+1} := down_j$.
         Increase $j$ by one.
\end{enumerate}
The iteration is stopped when $L(interm_k)$ nearly equals 
	 $A$ and $up_k$ and $down_k$ nearly equal $interm_k$
	 according to the adopted precision of approximation, and
	 then we set $t_0=interm_k$.

Having obtained $D_n$ we repeat the procedure above to find $u_n, D_{n+1}$ etc.
It is easily seen that $\lambda_n \leq \lambda_{n-1}.$
We iterate 
 until $| \lambda_n - \lambda_{n-1} | < \eps$, where $\eps$ is given.
   In the numerical experiments that we have done, we have taken 
   $\eps$ between $10^{-7} $ and $10^{-10}$.  

By the monotonicity of $\{\lambda_n\}$, the limit $\lim_{n\to\infty}
=\lambda_\infty$ exists.
However, it is not clear a priori whether $\lambda_\infty=
\Lambda_{\Omega}(\alpha,A)$ or not. In order to avoid the latter case,
we have repeated the same procedure with several different initial
shapes $D_0$.

The results of some of the computations that we have done are shown
in Figures 1-3.
They illustrate well Theorems
\ref{thtub}, \ref{thsymmpres} \ref{thsymmann}, \ref{thdumb},
and \ref{thasmallcxfb}.

\section{Some open problems and conjectures} \label{secopen}

In this section $D=D_{\alpha,A}$ will always denote an optimal configuration.

\begin{conj}(Uniqueness and convexity)
If $\Omega$ is convex then $D$ is unique,
and $D^c$ is convex (at least for $\alpha\leq\abar_\Omega(A)$).
\end{conj}
Concerning the restriction on $\alpha$ compare the remark
after Theorem \ref{thtub}. We have proved convexity for small
$\alpha$ in  Theorem \ref{thasmallcxfb}.

\begin{prob} (Regularity of the free boundary)
When is the boundary of an optimal configuration  smooth everywhere?
 In general, how can we control the size of 
singular sets of the free boundary? 
\end{prob}
In the convex case we have proved smoothness for small $\alpha$ in Theorem
\ref{thasmallcxfb}. A similar method should easily yield smoothness
of the free boundary for small $A$ and smooth $\partial\Omega$.

\begin{conj}
In dimension two the free boundary $\partial D$ is smooth outside
a finite set.
\end{conj}
We prove some restrictions on the singular set of $\partial D$
in [CGK].

\begin{prob}(Topology of $D$ and $D^c$)
If $\Omega$ is simply connected, is 
$D$  also connected 
even in the case $\alpha > \abar_\Omega(A)$ (cf. Theorem \ref{thtub})? 
If $A$ or $|\Omega|-A$ is small enough (with $\alpha$ fixed),
is $D^c$ always connected?
\end{prob}
Compare Proposition \ref{propAlarge} for the case of $A$ close to $|\Omega|$.
In a dumbbell $\psi_\Omega$ has two maxima. But numerical evidence
suggests the following conjecture:

\begin{conj} (One component of $D^c$ for dumbbell)
Let $\Omega_h$ be a dumbbell. Then for every $\alpha>0$ there is $\rho_0(\alpha,h)>0$
such that $D^c$ consists of one component (near one of the maxima of $\psi_\Omega$)
whenever $|\Omega|-A < \rho_0(\alpha,h)$. 
\end{conj}
Clearly, one would expect $\rho_0(\alpha,h)\to0$ as $\alpha\to 0$.

We now turn to questions of symmetry.  A very general problem is the following:
\begin{prob} (Symmetry and symmetry breaking)
Determine (at least qualitatively) the region in the space of parameters
where symmetry breaking occurs.

For annuli the parameters are $\alpha$, $\delta=A/|\Omega|$ and the ratio
$\tau$ of outer and inner radius (\lq thickness'). 
For dumbbells the parameters are $\alpha$, $A$
and the thickness of the handle.
\end{prob}
First results on this general problem are given by
Theorems \ref{thsymmann} and \ref{thdumb}.
The next three
conjectures address other aspects of this problem, i.e. they concern
other regions in parameter space. They are motivated by numerical
experiments.

\begin{conj} (Symmetry on dumbbells)
Let $\Omega_h$ be a dumbbell. Then for every $\alpha>0$ there is 
$\rho_1(\alpha,h)>0$ such that symmetry breaking occurs if and only if
$|\Omega|-A < \rho_1(\alpha,h)$.
\end{conj}

\begin{conj} (Symmetry on annuli)
For each $\alpha,\delta>0$ there is $\tau_0(\alpha,\delta)$ such that 
symmetry breaking occurs for the annulus of thickness $\tau$ if and only
if $\tau<\tau_0(\alpha,\delta)$.
\end{conj}
Theorem \ref{thsymmann} gives one half of this. The other half means that
the optimal configuration is rotationally symmetric for \lq thick' annuli.
Some aspects of this conjecture are discussed in [CGK].

More generally, it would be interesting to prove symmetry preservation in
{\em any} situation not covered by Theorem \ref{thsymmpres} (i.e. in a 
non-convex situation). In particular, a natural conjecture is:
\begin{conj}(Symmetry preservation for small $\alpha$)
For any domain $\Omega$ and any $A$ there is $\alpha_0(A,\Omega)$ such that for 
$\alpha\leq\alpha_0(A,\Omega)$ any optimal configuration
$D$ has the same symmetries as $\Omega$.
\end{conj}

Also, the analysis of the transition between the symmetric and asymmetric
situations would be interesting, as well as the shape of asymmetric
solutions for the annulus.

\begin{prob} (Relation between $D$ and the curvature of $\partial \Omega$)
Prove that $D$ is fat near points where $\partial \Omega$
has large positive curvature.
\end{prob}
For example see Figure 1.
For $\alpha=0$ and $A$ near zero this should
be not too hard. See [K1] for the case $\alpha=0$ under additional 
geometric assumptions.
From this one should obtain the result at least for small $\alpha$ and $A$
by perturbation.
In [CGK], Thm. 9, we prove in a model case that $D$ is thin near a portion
of the boundary which has large negative curvature.

\begin{prob}(Limit $\alpha\to\infty$)
Consider the restricted minimization problem, allowing only such
sets $D$ for which $D^c$ is  a  ball.
How does this relate to the limit $\alpha\to\infty$ in our problem?
Where does the center of an optimal ball lie?
\end{prob}
This is motivated as follows:
Formally, for $\alpha=\infty$ the eigenvalue
$\lambda_\Omega(\alpha,D)$ equals the first Dirichlet eigenvalue
of $D^c$. (The convergence to this value as $\alpha\to\infty$ is
proved in [HH] and [DKM], for example.)
Now by the Faber-Krahn inequality (see [Ch], for example), 
the first Dirichlet eigenvalue
of a domain of prescribed area is minimal if the domain is a ball.
So the optimal configuration for $\alpha$ large should be close to
a ball, at least when $A$ is close enough to $|\Omega|$ (so that
a ball of volume $|\Omega|-A$ fits into $\Omega$).

\begin{prob} (Other Elliptic Operators)
Consider the same optimization problem for a magnetic Schr\"{o}dinger operator 
$(i\nabla -\alpha\chi_D A(x))^2$
with constant magnetic field 
or 
a uniformly elliptic operator of divergence type
$-\nabla\{ (1+\alpha\chi_D(x))\nabla \}$. 
\end{prob}
We have no results for these operators, even if $\Omega$ is a ball. 
\section*{Appendix: Basic PDE facts} \label{secpdefacts}

Here we collect some well-known facts about uniform estimates for solutions 
of elliptic equations. We will state these for an equation
$$Pu=0,\quad P = \Delta + \sum_{j=1}^n b_j(x) \frac\partial{\partial x_j}
  + c(x),\quad x\in G, $$
where $P$ has measurable, uniformly bounded coefficients,
$u\in C^1(G)\cap C^0(\overline{G})$, and $G\subset\R^n$ is a
bounded open set.
In the following estimates, saying that the constants depend on $P$
will mean that they depend on $\sup_G (b_1,\ldots,b_n,c)$ and stay
bounded when this quantity stays bounded.

First, we have the uniform bound (see [GT, Thm.\ 8.15 and (8.38)])
\begin{equation} \label{eqpdeunif}
\sup_G |u| \leq C_{G,P} (\|u\|_{L^2(G)} + \sup_{\partial G} |u|).
\end{equation}

Second, we have Harnack's inequality: If $u\geq 0$ on $G$ and
$G'$ is a compact subset of $G$ then
\begin{equation} \label{eqpdehar}
\sup_{G'} u \leq C_{G,G',P} \inf_{G'} u.
\end{equation}

Combining these two we get the following slightly less standard estimate.
For $\eps>0$ let $G_\eps = \{x\in G:\, \dist(x,\partial G) \geq \eps\}.$

\begin{lem} \label{lhar}
For any $\eps>0$ there is a positive constant $c_{G,P,\eps}$ such that
for any $u\in C^1(G)\cap C^0(\overline{G})$ that solves $ Pu = 0$ 
and satisfies $u\geq 0$ one has
$$\inf_{G_{\eps}} u 
\geq c_{G,P,\eps} (\|u\|_{L^2(G)} - \sup_{\partial G} u). $$
\end{lem}
Here we set ${\displaystyle\inf_\emptyset u :=\infty}$.

\proof
We have
\begin{eqnarray*}
\int_G u^2 &=& \int_{G_\eps} u^2 + \int_{G\setminus G_\eps} u^2 
   \leq |G_\eps|\, \sup_{G_\eps} u^2 + |G\setminus G_\eps|\, \sup_G u^2 \\
  & \leq & C_{G,P,\eps} \inf_{G_\eps} u^2 + |G\setminus G_\eps|\,
                                        C'_{G,P} (\int_G u^2 + 
                                        \sup_{\partial G} u^2)
\end{eqnarray*}
where we used Harnack's inequality and the uniform estimate \eqref{eqpdeunif}.
If $\eps$ is so small that $|G\setminus G_{\eps}|\, C'_{G,P} < 1/2$
then we can subtract the last two terms, and the claim follows easily.
The claim for larger $\eps$ then follows from the fact that
$\inf_{G_{\eps'}} u \geq \inf_{G_\eps} u$ if $\eps'\geq\eps$.
\qed

\section*{REFERENCES}
{\footnotesize
[AB] Ashbaugh, M.S., Benguria, R.D., 
A sharp bound for the ratio of the first two 
eigenvalues of Dirichlet Laplacians and extensions, Ann. of Math., 
135(1992), 601--628.

[AH]  Ashbaugh, M.S., Harrell, E.,
Maximal and minimal eigenvalues and their associated nonlinear equations, J.Math.Phys., 28(1987), 1770-1786.

[AHS] Ashbaugh, M.S., Harrell, E., Svirsky, R., 
On minimal and maximal eigenvalue gaps 
and their causes, 
Pacific J.Math., 147(1991), 1--24.

[BL] Brascamp, H.J., Lieb, E., 
On extensions of the Brunn-Minkowski and Prekopa-Leindler 
theorems, including inequalities for log concave functions, and with an application to a diffusion equations, J.Fun.Ana., 22(1976), 366--389.

[BZ] Brothers, J.E., Ziemer, W.P., 
Minimal rearrangements of Sobolev functions, 
J. reine angew. Math., 384(1988), 153--179.

[CS] Caffarelli, L.A., Spruck, J., 
Convexity properties of solutions to some 
classical variational problems, 
Comm. in Partial Diff. Equ., 7(1982), 1337--1379.

[CGK] Chanillo, S., Grieser, D., and Kurata, K.,
The free boundary problem in the optimization of composite membranes.
Preprint.

[Ch] Chavel, I., {\it Eigenvalues in Riemannian Geometry}, Academic Press, 
New York, 1984. 

[C] Cox, S.J., The two phase drum with the deepest bass note,
Japan J. Indust. Appl. Math. 8(1991), 345-355.

[CL] Cox, S.J., Lipton, R., Extremal eigenvalue problems for two-phase
conductors, Arch. Rat. Mech. Anal. 136(1996), 101-117.

[CM] Cox, S.J., McLaughlin, J.R., Extremal eigenvalue problems for composite 
membranes, I and II, Appl. Math. Optim. 22(1990), 153--167, 169--187.
  
[DKM] Demuth, M., Kirsch, W., McGillivray, I., 
Schr\"{o}dinger operators-Geometric estimates in terms of the 
occupation time,
Comm. in Partial Diff. Equ., 20(1995), 37--57.

[Eg] Egnell, H., Extremal properties of the first eigenvalue 
of a class of elliptic eigenvalue problems, 
Ann. Scuol.Norm.Sup.Pisa, (1987), 1--48.

[FG] de Figueiredo, D.G., Gossez, J-P., 
Strict monotonicity of eigenvalues and unique continuation, Comm. P.D.E., 
17(1\& 2), (1992), 339--346.
 
[F] Friedman, A.,  On the regularity of the solutions of non-linear
elliptic and parabolic systems of partial differential equations,
J. Math. Mech. (now Indiana Math. J.) 7(1958), 43-60.

[GT] Gilbarg, D., Trudinger, N.T., 
{\it Elliptic Partial Differential Equations of Second Order}, Springer-Verlag, New York/Berlin, 1983.

[G] Giraud, G., Ann. Scient. de l'Ec. Norm. 43(1926), 1-128.

[HKK] Harrell, E., Kr\"{o}ger, P., Kurata, K., 
On the placement of an obstacle or a well so as to optimize the fundamental eigenvalue, preprint.

[HH] Hempel, P., Herbst, I., 
Strong magnetic fields, Dirichlet boundaries, and spectral gaps, Comm. Math. Phys., 169(1995), 237--260. 

[H] Hopf, E., \"Uber den funktionalen, insbesondere den analytischen
Character der L\"osungen elliptischer Differentialgleichungen zweiter
Ordnung, Math. Z. 34(1931), 194-233.

[K1] Kawohl. B., On the location of maxima of the gradient for
solutions to quasi-linear elliptic problems and a problem
raised by Saint Venant, J. of Elasticity 17(1987), 195-206.

[K2] Kawohl, B., Symmetrization -- or how to prove symmetry of
solutions to a PDE, in Partial Differential Equations, Theory and
Numerical Solution (W.J\"ager, J.Neusetal Eds.), Chapman \& Hall Res.
Notes in Math. 406(1999), 214-229.

[KL] Korevaar, N.J., Lewis, J.L., Convex solutions of certain elliptic
equations have constant rank Hessians, Arch. Rat. Mech. Anal. 97(1987),
19-32.

[Kr] Krein, M.G., On certain problems on the maximum and minimum 
of characteristic values and on the Lyapunov zones of stability, AMS Translations 
Ser. 2, 1(1955), 163--187.

[LL] Lieb, E., Loss, M., {\it Analysis}, Amer. Math. Soc., 1997. 

[M] Morrey, C.B., On the analyticity of solutions of analytic
non-linear elliptic systems of PDE I, Amer. J. Math. 80(1958), 198-218.

[Pi] Pironneau, O., 
{\it Optimal shape design for elliptic systems}, Springer-Verlag, New York Inc., 
1984.

[TK] Tsukuda, Y.,  Kaizu, S.,  
Proceedings of the tenth symposium of numerical fluid mechanics, (1996), 
220--221 (in Japanese with English abstract).

} 


\end{document}